\documentclass[12]{article}
\usepackage{amssymb}	
\usepackage{amsmath}    
\usepackage{mathtools}
\usepackage{mhsetup}

\usepackage{setspace} 
\usepackage{latexsym}	
\usepackage{longtable}
\usepackage{multirow}
\usepackage{epsfig}

\usepackage{wrapfig}
\usepackage{graphicx}

\allowdisplaybreaks

\newtheorem{algorithm}{Algorithm}[section]

\newtheorem{remark}{Remark}[section]

\newcommand{\qed}{\nobreak \ifvmode \relax \else \ifdim\lastskip<1.5em \hskip-\lastskip \hskip1.5em plus0em minus0.5em \fi \nobreak \vrule height0.75em width0.5em depth0.25em\fi} 

\def\A{{\bf A}}
\def\B{{\bf B}}

\def\E{{\bf E}}
\def\F{{\bf F}}

\def\I{{\bf I}}
\def\J{{\bf J}}
\def\K{{\bf K}}

\def\0{{\bf 0}}
\def\P{{\bf P}}
\def\Q{{\bf Q}}
\def\R{{\bf R}}
\def\S{{\bf S}}
\def\T{{\bf T}}

\def\W{{\bf W}}
\def\X{{\bf X}}

\def\b{{\bf b}}

\def\d{{\bf d}}
\def\e{{\bf e}}
\def\f{{\bf f}}
\def\g{{\bf g}}
\def\h{{\bf h}}
\def\i{{\bf i}}

\def\m{{\bf m}}

\def\q{{\bf q}}

\def\t{{\bf t}}
\def\u{{\bf u}}
\def\v{{\bf v}}

\def\x{{\bf x}}

\def\Tr{{\rm T}}

\def\diag{{\rm diag}}

\title{Spacecraft Attitude and Reaction Wheel Desaturation Combined Control Method}
\author {Yaguang Yang\thanks{
Office of Research, NRC, 21 Church Street, Rockville, 20850. Email:
yaguang.yang@verizon.net} 
}
\date{\today}

\begin{document}

\maketitle  

\begin{abstract}
Two popular types of spacecraft actuators are reaction wheels and 
magnetic torque coils. Magnetic torque coils are particularly 
interesting because they can be used for both attitude control 
and reaction wheel momentum management (desaturation control). 
Although these two tasks
are performed at the same time using the same set of actuators, most design methods deal with only one of the these tasks or 
consider these two tasks separately. In this paper, a design with 
these two tasks in mind is formulated as a single problem. A 
periodic time-varying linear quadratic regulator design 
method is then proposed to solve this problem. A simulation example is provided to describe the benefit of the new strategy.
\end{abstract}

{\bf Keywords:} Spacecraft attitude control, reaction wheel 
desaturation, linear time-varying system, reduced quaternion 
model, linear quadratic regulator.

\newpage

\section{ Introduction}

Spacecraft attitude control and reaction wheel desaturation are 
normally regarded as two different control system design 
problems and are discussed in separate chapters in text books, 
such as \cite{wertz78,sidi97}. While spacecraft attitude control 
using magnetic torques has been one of the main research areas 
(see, for example, \cite{sl05,rmb15} and extensive
references therein), there are quite a few research papers that 
address reaction wheel momentum management, see for example, 
\cite{dbp90,ch99,gqt06} and references therein. In \cite{dbp90},
Dzielsk et al. formulated the problem as an optimization
problem and a nonlinear programming method was proposed to 
find the solution. Their method can be very expensive and
there is no guarantee to find the global optimal solution.
Chen et al. \cite{ch99} discussed optimal desaturation 
controllers using magnetic torques and thrusters. Their methods 
find the optimal torques which, however, may not be able to 
achieve by magnetic torque coils because given the desired torques
in a three dimensional space, magnetic torque coils can only generate 
torques in a two  dimensional plane \cite{sidi97}. 
Like most publications on this problem, the above two papers do
not consider the time-varying effect of the geomagnetic field
in body frame, which arises when a spacecraft flies around the 
Earth. Giulietti et al. \cite{gqt06} considered the same problem 
with more details on the geomagnetic field, but the periodic 
feature of the magnetic field along the orbit was not used in 
their proposed design.
In addition, all these proposed designs considered only momentum
management but not attitude control.

Since both attitude control and reaction wheel desaturation are 
performed at the same time using the same magnetic torque coils,
the control system design should consider these two design 
objectives at the same time and some very recent research papers tackled the problem in this direction, for example, 
\cite{tappz15,agra15}. In \cite{tappz15}, Tregouet et al. studied the problem of the spacecraft stabilization and reaction wheel 
desaturation at the same time. They considered time-variation of
the magnetic field in body frame, and their reference frame was 
the inertial frame. However, for a Low Earth Orbit (LEO) 
spacecraft that uses Earth's magnetic field, the reference frame for the spacecraft is most likely Local Vertical Local Horizontal 
(LVLH) frame. In addition, their design method depends on some 
assumption which is not easy to verify and their 
proposed design does not use the periodic feature
of the magnetic field. Moreover, their design is
composed of two loops, which is essentially an idea of dealing
with attitude control and wheel momentum management in separate
considerations. In \cite{agra15}, a heuristic proportional
controller was proposed and a Lyapunov function was used to 
prove that the controller can simultaneously stabilize the 
spacecraft with respect to the LVLH frame and achieve reaction 
wheel management. But this design method does not consider the
the time-varying effect of the geomagnetic field in body frame.
Although these two designs are impressive, as we have seen, 
these designs do not consider some factors in reality and their 
solutions are not optimal.

In this paper, we propose a more attractive design method
which considers as many factors as practical. The controlled 
attitude is aligned with LVLH frame. A general reduced 
quaternion model, including (a) reaction wheels, (b) magnetic 
torque coils, (c) the gravity gradient torque, and (d) the periodic 
time-varying effects of the geomagnetic field along the orbit
and its interaction with magnetic torque coils, is proposed. 
The model is an extension of the one discussed in 
\cite{yang10}. A single objective function, which considers 
the performance of both attitude control and 
reaction wheel management at the same time, is suggested. Since
a well-designed periodic controller for a period system is better 
than constant controllers as pointed out in \cite{flamm91,kpt85},
this objective function is optimized using the solution of
a matrix periodic Riccati equation described
in \cite{yang16a}, which leads to a periodic time-varying 
optimal control. Based on the algorithm for the periodic
Riccati equations \cite{yang16a}, we show that the design can 
be calculated in an efficient way and the designed 
controller is optimal for both the spacecraft attitude control 
and for the reaction wheel momentum manage at the same time. 
We provide a simulation test to demonstrate that the designed 
system achieves more accurate attitude than the optimal
control system that uses only magnetic torques. Moreover, the
designed controller based on LQR method works on the nonlinear
spacecraft system.

The remainder of the paper is organized as follows. Section 2 
derives the reduced quaternion spacecraft control system model using reaction
wheel and magnetic control torques with the attitude defined as 
the rotation of the body frame respect to the LVLH frame. 
Section 3 reduces the nonlinear spacecraft system model to a 
linearized periodic time-varying model which includes the
time-varying geomagnetic field along the orbit, the gravity 
gradient disturbance torque, the reaction wheel speed control, 
and the magnetic torque control. Section 4 introduces a
single objective function for both attitude control 
and wheel management. It also gives the optimal control solutions
for this linear time-varying system in different conditions. Simulation test is provided
in Section 5. The conclusions are summarized in Section 6.

\section{Spacecraft model for attitude and reaction wheel desaturation control} 

Throughout the discussion, we assume that the inertia matrix of a 
spacecraft $\J=\diag(J_1, J_2, J_3)$ is a diagonal 
matrix. This assumption is reasonable because in practical
spacecraft design, spacecraft inertia matrix J is always designed as close to a diagonal matrix as possible \cite{yang12a}.
(It is actually very close to a diagonal matrix.)
For spacecraft using Earth's magnetic torques, the nadir pointing model
is probably the mostly desired one. Therefore, the attitude of the spacecraft is represented by the 
rotation of the spacecraft body frame relative to the local vertical and local horizontal 
frame. This means that the quaternion and spacecraft body 
rate should be represented in terms of the rotation of the 
spacecraft body frame relative to the LVLH frame (see 
\cite{yang12a} for the definition of LVLH frame). 

Let $\boldsymbol{\omega}=[\omega_1,\omega_2, \omega_3]^{\Tr}$ be
the body rate with respect to the LVLH frame represented in the 
body frame, $\boldsymbol{\omega}_{lvlh}=[0,\omega_0, 0]^{\Tr}$ 
the orbit rate (the rotation of LVLH frame) with respect 
to the inertial frame represented in the LVLH frame\footnote{
For a circular orbit, given the spacecraft orbital period around the 
Earth $P$, $\omega_0 = \frac{2\pi}{P}$ is a known constant.}, and 
$\boldsymbol{\omega}_I=[\omega_{I1},\omega_{I2}, \omega_{I3}]^{\Tr}$ 
be the angular velocity vector of the spacecraft
body with respect to the inertial frame, represented in the 
spacecraft body frame. Let $A_l^b$ represent the rotational
transformation matrix from the LVLH frame to
the spacecraft body frame. Then, $\boldsymbol{\omega}_I$ can be expressed as 
\cite{yang10,yang12a}
\begin{equation}
\boldsymbol{\omega}_{I}=\boldsymbol{\omega} + \A_{l}^{b} 
\boldsymbol{\omega}_{lvlh} = \boldsymbol{\omega} + 
\boldsymbol{\omega}_{lvlh}^b,
\label{3}
\end{equation}
where $\boldsymbol{\omega}_{lvlh}^b$ is the rotational rate of
LVLH frame relative to the inertial frame represented in the
spacecraft body frame. Assuming that the orbit is circular,
i.e., $\dot{\boldsymbol{\omega}}_{lvlh}=0$,
using the fact (see \cite[eq.(19)]{yang12a}) 
\begin{equation}
\dot{\A}_{l}^{b} = -\boldsymbol{\omega} \times {\A}_{l}^{b},
\end{equation}
we have
\begin{eqnarray}
\dot{\boldsymbol{\omega}}_I &=& \dot{\boldsymbol{\omega}}
 + \dot{\A}_{l}^{b} \boldsymbol{\omega}_{lvlh} + \A_{l}^{b} 
\dot{\boldsymbol{\omega}}_{lvlh}
\nonumber \\
&= & \dot{\boldsymbol{\omega}} -\boldsymbol{\omega}
 \times {\A}_{l}^{b} \boldsymbol{\omega}_{lvlh} = 
\dot{\boldsymbol{\omega}} -\boldsymbol{\omega}
 \times\boldsymbol{\omega}_{lvlh}^b.
\label{4}
\end{eqnarray}

Assuming that the three reaction wheels are aligned with the body frame axes, 
the total angular momentum of the spacecraft $\h_T$ in the body 
frame comprises the angular momentum of the 
spacecraft $\J \boldsymbol{\omega}_I$ and the angular momentum of the reaction wheels $\h_w=[h_{w1},h_{w2},h_{w3}]^{\Tr}$ 
\begin{equation}
\h_T =\J \boldsymbol{\omega}_I + \h_w,
\label{angularM}
\end{equation}
where 
\begin{equation}
{\h}_w= \J_w  {\boldsymbol{\Omega}},
\end{equation} 
$\J_w =\diag(\J_{w_1},\J_{w_2},\J_{w_3})$ is the 
inertia matrix of the three reaction wheels aligned with the
spacecraft body axes, and 
$\boldsymbol{\Omega} = [\Omega_1, \Omega_2,\Omega_3]^{\Tr}$
is the angular rate vector of the three reaction wheels. Let $\h_T'$ be the same 
vector of $\h_T$ represented in inertial frame. Let $\t_T$ 
be the total external torques acting on the spacecraft, 
we have (see \cite{serway04})
$
\t_T = \frac{d \h_T'}{dt} \bigg|_b.
\label{basic}
$
Using eq. (20) of \cite{yang12a} and equation (\ref{angularM}), 
we have the dynamics equations of the spacecraft as follows
\begin{eqnarray}
\J \dot{\boldsymbol{\omega}}_I  + \dot{\h}_w & = &
\left( \frac{d \h_T}{dt} \right)\bigg|_b =
-\boldsymbol{\omega}_I  \times \h_T+\left( \frac{d \h_T'}{dt} \right)\bigg|_b
\nonumber \\
& = &
- \boldsymbol{\omega}_I  \times (\J \boldsymbol{\omega}_I 
 + \h_w) + \t_T,
\label{dynamics}
\end{eqnarray}
where $\t_T$ includes the gravity gradient torque $\t_g$, magnetic 
control torque $\t_m$, and internal and external disturbance torque $\t_d$ (including residual magnetic moment induced torque, 
atmosphere induced torque, solar radiation torque, etc). The 
torques generated by the reaction wheels $\t_w$ are given by 
\[
\t_w= \dot{\h}_w= \J_w  \dot{\boldsymbol{\Omega}}.
\]
Substituting these relations into (\ref{dynamics}) gives
\begin{equation}
\J \dot{\boldsymbol{\omega}}_I  =
- \boldsymbol{\omega}_I  \times (\J \boldsymbol{\omega}_I  + 
\J_w  {\boldsymbol{\Omega}} ) -\t_w + \t_g+\t_m+\t_d.
\label{dynamics1}
\end{equation}
Substituting (\ref{3}) and (\ref{4}) into (\ref{dynamics1}), we 
have
\begin{equation}
\J \dot{\boldsymbol{\omega}} = \J \boldsymbol{\omega}
 \times\boldsymbol{\omega}_{lvlh}^b 
- (\boldsymbol{\omega} + \boldsymbol{\omega}_{lvlh}^b )
 \times [ \J ( \boldsymbol{\omega} + \boldsymbol{\omega}_{lvlh}^b)
+ \J_w  {\boldsymbol{\Omega}} ] -\t_w + \t_g+\t_m+\t_d.
\label{dynamics2}
\end{equation}
Let 
\begin{equation}
\bar{\q}=[q_0, q_1, q_2, q_3]^{\Tr}=[q_0, \q^{\Tr}]^{\Tr}=
\left[ \cos(\frac{\alpha}{2}), \hat{\e}^{\Tr}\sin(\frac{\alpha}{2}) \right]^{\Tr}
\end{equation}
be the quaternion representing the rotation of the body frame relative to the LVLH frame, 
where $\hat{\e}$ is the unit length rotational axis and $\alpha$ is the rotation angle about $\hat{\e}$. 
Therefore, the reduced kinematics equation becomes \cite{yang10}
\begin{eqnarray} \nonumber
\left[  \begin{array} {c} \dot{q}_1 \\ \dot{q}_2 \\ \dot{q}_3
\end{array} \right]
& = & \frac{1}{2}  \left[  \begin{array} {ccc} 
\sqrt{1-q_1^2-q_2^2-q_3^2} & -q_3  & q_2 \\
q_3 & \sqrt{1-q_1^2-q_2^2-q_3^2} & -q_1 \\
-q_2  &  q_1  & \sqrt{1-q_1^2-q_2^2-q_3^2} \\
\end{array} \right] 
\left[  \begin{array} {c}  \omega_{1} \\ \omega_{2} \\ \omega_{3}
\end{array} \right]  \\
& = & \g(q_1,q_2, q_3, \boldsymbol{\omega}),
\label{nadirModel2}
\end{eqnarray}
or simply
\begin{equation}
\dot{\q}=\g(\q, \boldsymbol{\omega}).
\label{gfunction}
\end{equation}
Since (see \cite{yang10,yang12a}), 
\begin{eqnarray}  \nonumber
\A_l^b = \left[ \begin{array} {ccc}
2 q_0^2-1+2q_1^2  &  2q_1q_2+2q_0q_3 & 2q_1q_3-2q_0q_2 \\
2q_1q_2-2q_0q_3   & 2q_0^2-1+2q_2^2  & 2q_2q_3+2q_0q_1 \\
2q_1q_3+2q_0q_2   &  2q_2q_3-2q_0q_1 & 2q_0^2-1+2q_3^2  
\end{array} \right],
\end{eqnarray}
we have
\begin{equation}
\boldsymbol{\omega}_{lvlh}^b=\A_l^b\boldsymbol{\omega}_{lvlh}
=\left[ \begin{array}{c}
2q_1q_2+2q_0q_3 \\ 2q_0^2-1+2q_2^2 \\ 2q_2q_3-2q_0q_1
\end{array} \right] \omega_0,
\label{n9}
\end{equation}
which is a function of $\q$. Interestingly, given spacecraft 
inertia matrix $\J$, $\t_g$ is also a function of $\q$. Using the
facts (a) the spacecraft mass is negligible compared to the 
Earth mass, and (b) the size of the spacecraft is negligible 
compared to the magnitude of the vector from the center of the 
Earth to the center of the mass of the spacecraft $\R$, the 
gravitational torque is given by \cite[page 367]{wie98}:
\begin{align}
\t_g = \frac{3 \mu }{|\R|^{5}} \R  \times \J \R,
\label{gg1}
\end{align}
where $\mu = G M$, $G=6.669*10^{-11}m^3/kg-s^2$ is the universal 
constant of gravitation, and $M$ is the mass of the Earth. 
Noticing that in local vertical local horizontal frame, 
$\R_l = \left[ 0,0,-|\R| \right]^{\Tr}$, we can represent $\R$ in body frame as 
\begin{equation}
\R = \A_l^b \R_l =
\left[ \begin{array} {ccc}
2 q_0^2-1+2q_1^2  &  2q_1q_2+2q_0q_3 & 2q_1q_3-2q_0q_2 \\
2q_1q_2-2q_0q_3   & 2q_0^2-1+2q_2^2  & 2q_2q_3+2q_0q_1 \\
2q_1q_3+2q_0q_2   &  2q_2q_3-2q_0q_1 & 2q_0^2-1+2q_3^2 \end{array} \right]
\left[ \begin{array}{c} 0 \\ 0 \\ -|\R| \end{array} \right].
\label{rl}
\end{equation}
Denote the last column of $\A_l^b$ as $ \A_l^b(:,3) $, and using the following relation \cite[page 109]{sidi97}
\begin{equation}
\omega_0 =\sqrt{\frac{\mu}{|\R|^3}}
\label{aOrbitV}
\end{equation}
and (\ref{rl}), we can rewrite (\ref{gg1}) as
\begin{align}
\t_g = 3 \omega_0^2 \A_l^b(:,3)   \times \J \A_l^b(:,3).
\label{ginterM1}
\end{align}

Let $\b(t)=[b_1(t),b_2(t),b_3(t)]^{\Tr}$ be the Earth's magnetic 
field in the spacecraft coordinates, computed using the 
spacecraft position, the spacecraft attitude, and a spherical harmonic model of the Earth's magnetic field \cite{wertz78}. Let 
$\m=[m_1,m_2,m_3]^{\Tr}$ be the spacecraft magnetic torque coils'
induced magnetic moment in the spacecraft coordinates. The desired  magnetic control torque $\t_m$  may not be achievable because 
\begin{equation}
\t_m = \m \times \b = - \b \times \m
\end{equation}
provides only a torque in a two dimensional plane but not in the 
three dimensional space \cite{sidi97}. However, the spacecraft 
magnetic torque  coils' induced magnetic moment $\m$ is an 
achievable engineering variable. Therefore, equation 
(\ref{dynamics2}) should be rewritten as 
\begin{equation}
\J \dot{\boldsymbol{\omega}} = \f(\boldsymbol{\omega},\boldsymbol{\Omega},\q) -\t_w + \t_g - \b \times \m +\t_d,
\label{dynamics3}
\end{equation}
where 
\begin{equation}
\f(\boldsymbol{\omega},\boldsymbol{\Omega},\q)
=\J \boldsymbol{\omega}
 \times\boldsymbol{\omega}_{lvlh}^b 
- (\boldsymbol{\omega} + \boldsymbol{\omega}_{lvlh}^b )
 \times [ \J ( \boldsymbol{\omega} + \boldsymbol{\omega}_{lvlh}^b)
+ \J_w {\boldsymbol{\Omega}} ].
\label{ffunction}
\end{equation}
Notice that the cross product of $\b \times \m$ can be expressed 
as product of an asymmetric matrix $\b^{\times}$ and the vector $\m$ with
\begin{equation}
\b^{\times}= \left[ \begin{array}{ccc}
0 & -b_3 & b_2 \\
b_3 & 0 & -b_1 \\ 
-b_2  & b_1 & 0
\end{array} \right].
\label{productM}
\end{equation}
Denote the system states $\x=[\boldsymbol{\omega}^{\Tr},\boldsymbol{\Omega}^{\Tr},\q^{\Tr}]^{\Tr}$ and control inputs 
$\u=[\t_w^{\Tr},\m^{\Tr}]^{\Tr}$, the spacecraft control system 
model can be written as follows:
\begin{subequations}
\begin{gather}
\J \dot{\boldsymbol{\omega}} = \f(\boldsymbol{\omega},\boldsymbol{\Omega},\q) + \t_g-[\I, \b^{\times}] \u +\t_d, 
\label{nonlinearA} \\
\J_w  \dot{\boldsymbol{\Omega}} = \t_w, \label{nonlinearB}  \\
\dot{\q}=\g(\q, \boldsymbol{\omega}) \label{nonlinearC}.
\end{gather}
\label{dynamicsAll}
\end{subequations}

\begin{remark}
The reduced quaternion, instead of the full quaternion, is proposed in this model
because of many merits discussed in \cite{yang10,yang12,yang14}. 
\end{remark}

\section{Linearized model for attitude and reaction wheel desaturation control}

The nonlinear model of (\ref{dynamicsAll}) can be used to design 
control systems. One popular design method for nonlinear model 
involves Lyponuv stability theorem, which is actually used in \cite{tappz15,agra15}. A design based on this method focuses 
on stability but not on performance. Another widely known method 
is nonlinear optimal control design \cite{dbp90}, it normally 
produces an open loop controller which is not robust 
\cite{lvs12} and its computational cost is high. Therefore,
We propose to use Linear Quadratic Regulator (LQR) which achieves
the optimal performance for the linearized system and is
a closed-loop feedback control. 
Our task in this section is to derive the linearized model for the nonlinear 
system (\ref{dynamicsAll}).

Using the linearization technique of \cite{yang10,yang12a},
we can express $\boldsymbol{\omega}_{lvlh}^b $ in (\ref{n9})
approximately as a linear function of $\q$ as follows
\begin{eqnarray}
\boldsymbol{\omega}_{lvlh}^b
\approx
\left[ \begin{array}{c}
2q_3 \\ 1 \\ -2q_1
\end{array} \right] \omega_0
= \left[  \begin{array}{ccc} 
0 & 0 & 2\omega_0 \\
0 & 0 & 0 \\
-2\omega_0 & 0 & 0 \\
\end{array}  \right] \q 
+ \left[  \begin{array}{c} 
0 \\
\omega_0\\
0 \\
\end{array}  \right]. 
\label{9}
\end{eqnarray}
Similarly, we can express $\t_g$ in (\ref{ginterM1}) approximately 
as a linear function of $\q$ as follows
\begin{eqnarray}
\t_g \approx \left[  \begin{array}{c} 
6\omega_0^2 (J_{3}-J_{2}) q_1 \\
6\omega_0^2 (J_{3}-J_{1}) q_2 \\
0 \\
\end{array}  \right]= \left[  \begin{array}{ccc} 
6\omega_0^2 (J_{3}-J_{2}) &  0 & 0 \\
0 & 6\omega_0^2 (J_{3}-J_{1})  & 0 \\
0 & 0 & 0 \\
\end{array}  \right] \q := \T \q.
\label{22}
\end{eqnarray}
Since $\t_g$ and $\boldsymbol{\omega}_{lvlh}^b$ are
functions of $\q$, the linearized spacecraft model can be expressed as follows:
\begin{eqnarray}
\left[ \begin{array}{ccc}
\J & \0 & \0 \\
\0 & \J_w  & \0 \\ 
\0  & \0 & \I
\end{array} \right]
\left[ \begin{array}{c}
\dot{\boldsymbol{\omega}}  \\
\dot{\boldsymbol{\Omega}} \\
\dot{\q}
\end{array} \right]
=  \left[ \begin{array}{ccc}
\frac{\partial{\f}}{\partial{\boldsymbol{\omega}}}  
& \frac{\partial{\f}}{\partial{\boldsymbol{\Omega}}}   
& \frac{\partial{\f}}{\partial{\boldsymbol{q}}}  + \T \\
\0 & \0  & \0 \\ 
\frac{\partial{\g}}{\partial{\boldsymbol{\omega}}}  & \0 & \frac{\partial{\g}}{\partial{\q}}
\end{array} \right]
\left[ \begin{array}{c}
{\boldsymbol{\omega}}  \\
{\boldsymbol{\Omega}} \\
{\q}
\end{array} \right]
+\left[ \begin{array}{cc}
-\I & -\b^{\times} \\
\I  & \0 \\ 
\0  & \0 
\end{array} \right]
\left[ \begin{array}{c}
\t_w  \\
\m
\end{array} \right]
+ \left[ \begin{array}{c}
\t_d  \\
\0   \\
\0
\end{array} \right],
\label{linearModel}
\end{eqnarray}
where $\frac{\partial{\f}}{\partial{\boldsymbol{\omega}}}$, 
$\frac{\partial{\f}}{\partial{\boldsymbol{\Omega}}}$, 
$\frac{\partial{\f}}{\partial{\boldsymbol{q}}}$, $\frac{\partial{\g}}{\partial{\boldsymbol{\omega}}}$, 
and $\frac{\partial{\g}}{\partial{\q}}$ are evaluated at the
desired equilibrium point ${\boldsymbol{\omega}}=0$, 
${\boldsymbol{\Omega}} =0$, and ${\q}=0$.
Using the definition of (\ref{productM}), (\ref{9}), (\ref{22}), and (\ref{ffunction}), we have 
\begin{eqnarray}
\frac{\partial{\f}}{\partial{\boldsymbol{\omega}}}  
 \bigg|_{\substack{\boldsymbol{\omega} \approx 0 \\ \boldsymbol{\Omega} \approx 0 \\ \q \approx 0}} 
&  \approx  & -\J (\boldsymbol{\omega}_{lvlh}^b)^{\times}+  (\J\boldsymbol{\omega}_{lvlh}^b)^{\times}
 -  (\boldsymbol{\omega}_{lvlh}^b)^{\times} \J
 \bigg|_{\substack{\boldsymbol{\omega} \approx 0 \\ \boldsymbol{\Omega} \approx 0 \\ \q \approx 0}}
 \nonumber \\
 & = & -\J \left[ \begin{array}{ccc}
0 & 0 & \omega_0 \\
0 & 0 & 0 \\ 
-\omega_0  & 0 & 0
\end{array} \right]
+ \left[ \begin{array}{ccc}
0 & 0 & J_2 \omega_0 \\
0 & 0 & 0 \\ 
-J_2 \omega_0  & 0 & 0
\end{array} \right]
- \left[ \begin{array}{ccc}
0 & 0 & \omega_0 \\
0 & 0 & 0 \\ 
-\omega_0  & 0 & 0
\end{array} \right] \J
 \nonumber \\
 & = &  \left[ \begin{array}{ccc}
0 & 0 & \omega_0(-J_1+J_2-J_3) \\
0 & 0 & 0 \\ 
\omega_0(J_1-J_2+J_3)  & 0 & 0
\end{array} \right],
\label{12}
\end{eqnarray}
\begin{eqnarray}
\frac{\partial{\f}}{\partial{\boldsymbol{\Omega}}}  
 \bigg|_{\substack{\boldsymbol{\omega} \approx 0 \\ \boldsymbol{\Omega} \approx 0 \\ \q \approx 0}}
&  \approx  & -(\boldsymbol{\omega})^{\times}\J_w -  (\boldsymbol{\omega}_{lvlh}^b)^{\times}\J_w
 \bigg|_{\substack{\boldsymbol{\omega} \approx 0 \\ \boldsymbol{\Omega} \approx 0 \\ \q \approx 0}}
= - \left[ \begin{array}{ccc}
0 & 0 & \omega_0 \\
0 & 0 & 0 \\ 
-\omega_0  & 0 & 0
\end{array} \right] \J_w
 \nonumber \\
 & = &  \left[ \begin{array}{ccc}
0 & 0 & -\omega_0 J_{w_3} \\
0 & 0 & 0 \\ 
\omega_0 J_{w_1}  & 0 & 0
\end{array} \right],
\label{12a}
\end{eqnarray}
and
\begin{eqnarray}
\frac{\partial{\f}}{\partial{\boldsymbol{q}}}  
 \bigg|_{\substack{\boldsymbol{\omega} \approx 0 \\ \boldsymbol{\Omega} \approx 0 \\ \q \approx 0}}
&  \approx  & - \frac{\partial}{\partial{\boldsymbol{q}}}  \left( \boldsymbol{\omega}_{lvlh}^b 
 \times \J  \boldsymbol{\omega}_{lvlh}^b  \right)
 \bigg|_{\substack{\boldsymbol{\omega} \approx 0 \\ \boldsymbol{\Omega} \approx 0 \\ \q \approx 0}}
 \nonumber \\
 &   \approx   &   (\J  \boldsymbol{\omega}_{lvlh}^b)^{\times}
 \left[ \begin{array}{ccc}
0 & 0 & 2 \omega_0 \\
0 & 0 & 0 \\ 
-2 \omega_0  & 0 & 0
\end{array} \right] 
- ( \boldsymbol{\omega}_{lvlh}^b )^{\times} \J  
 \left[ \begin{array}{ccc}
0 & 0 & 2 \omega_0 \\
0 & 0 & 0 \\ 
-2 \omega_0  & 0 & 0
\end{array} \right] 
 \nonumber \\
 &   \approx  &  \omega_0 \left(  \left[ \begin{array}{c}
2J_1q_3  \\
J_2 \\ 
-2J_3q_1
\end{array} \right]^{\times}
-\left[ \begin{array}{ccc}
0 & 0 & 1  \\
0 & 0 & 0 \\ 
-1   & 0 & 0
\end{array} \right] \J
\right)
\left[ \begin{array}{ccc}
0 & 0 & 2\omega_0  \\
0 & 0 & 0 \\ 
-2\omega_0   & 0 & 0
\end{array} \right]
 \nonumber \\
 &   \approx   &  \omega_0 \left[ \begin{array}{ccc}
0 & 0 & J_2-J_3  \\
0 & 0 & 0 \\ 
J_1-J_2   & 0 & 0
\end{array} \right] 
\left[ \begin{array}{ccc}
0 & 0 & 2\omega_0  \\
0 & 0 & 0 \\ 
-2\omega_0   & 0 & 0
\end{array} \right]
 \nonumber \\
 &  =   &   
\left[ \begin{array}{ccc}
2\omega_0^2  ( J_3-J_2)  & 0 & 0\\
0 & 0 & 0 \\ 
0  & 0 & 2\omega_0^2  ( J_1-J_2) 
\end{array} \right].
\label{12b}
\end{eqnarray}
From (\ref{gfunction}), we have
\begin{eqnarray}   
\frac{\partial{\g}}{\partial{\boldsymbol{\omega}}} 
\bigg|_{\substack{\boldsymbol{\omega} \approx 0 \\ \q \approx 0}} 
 \approx \frac{1}{2} \I,
\label{14}
\end{eqnarray}
\begin{eqnarray}   
\frac{\partial{\g}}{\partial{\q}} 
\bigg|_{\substack{\boldsymbol{\omega} \approx 0 \\ \q \approx 0}} 
 \approx   \0.
\label{15}
\end{eqnarray}
Substituting (\ref{22}), (\ref{productM}), (\ref{12}), 
(\ref{12a}), (\ref{12b}), (\ref{14}), and (\ref{15}) into
(\ref{linearModel}), we have
\begin{eqnarray}  
\left[ \begin{array}{c}
\dot{\boldsymbol{\omega}}  \\
\dot{\boldsymbol{\Omega}} \\
\dot{\q}
\end{array} \right]
& =  &
 \left[ \begin{array}{ccc}
\J^{-1} \frac{\partial{\f}}{\partial{\boldsymbol{\omega}}}  
& \J^{-1}  \frac{\partial{\f}}{\partial{\boldsymbol{\Omega}}}   
& \J^{-1} \left( \frac{\partial{\f}}{\partial{\boldsymbol{q}}}  +  \T \right) \\
\0 & \0  & \0 \\ 
\frac{\partial{\g}}{\partial{\boldsymbol{\omega}}}  & \0 & \frac{\partial{\g}}{\partial{\q}}
\end{array} \right]
\left[ \begin{array}{c}
{\boldsymbol{\omega}}  \\
{\boldsymbol{\Omega}} \\
{\q}
\end{array} \right]
+\left[ \begin{array}{cc}
-\J^{-1}  & -\J^{-1} \b^{\times} \\
\J_w^{-1}  & \0 \\ 
\0  & \0 
\end{array} \right]
\left[ \begin{array}{c}
\t_w  \\
\m
\end{array} \right]
+ \left[ \begin{array}{c}
\J^{-1} \\
\0  \\
\0
\end{array} \right] \t_d 
\nonumber \\  
& = & \left[ \begin{array}{ccccccccc}
0  &  0  &  \omega_0\frac{J_1-J_2+J_3}{-J_1}   &  0  &  0  &  \frac{\omega_0 J_{w_3}}{-J_1}   &  8\omega_0^2 \frac{J_3-J_2}{J_1}  &  0  &  0 \\
0  &  0  &  0  &  0  &  0  &  0   &  0  &  6\omega_0^2 \frac{J_3-J_1}{J_2}  &  0 \\
 \omega_0\frac{J_1-J_2+J_3}{J_3}  &  0  &  0  &  \frac{\omega_0 J_{w_1}}{J_3}  &  0  &  0   &  0  &  0  &  2\omega_0^2 \frac{J_1-J_2}{J_3} \\
0  &  0  &  0  &  0  &  0  &  0   &  0  &  0  &  0 \\
0  &  0  &  0  &  0  &  0  &  0   &  0  &  0  &  0 \\
0  &  0  &  0  &  0  &  0  &  0   &  0  &  0  &  0 \\
0.5  &  0  &  0  &  0  &  0  &  0   &  0  &  0  &  0 \\
0  &  0.5  &  0  &  0  &  0  &  0   &  0  &  0  &  0 \\
0  &  0  &  0.5  &  0  &  0  &  0   &  0  &  0  &  0 
\end{array} \right] 
\left[ \begin{array}{c}
\omega_1 \\
\omega_2 \\
\omega_3 \\
\Omega_1 \\
\Omega_2 \\
\Omega_3 \\
q_1  \\
q_2  \\
q_3
\end{array} \right] 
 \nonumber \\
& +   &  \left[ \begin{array}{cccccc}
-J_1^{-1}  &  0  &  0  &  0  &  \frac{b_3}{J_1}  &  -\frac{b_2}{J_1}       \\
0  &  -J_2^{-1}   &  0  &  -\frac{b_3}{J_2}  &  0  &  \frac{b_1}{J_2}     \\
0  &  0  &  -J_3^{-1}   & \frac{b_2}{J_3}  & -\frac{b_1}{J_3}   &  0       \\
J_{w_1}^{-1}   &  0  &  0  &  0  &  0  &  0       \\
0  &  J_{w_2}^{-1}   &  0  &  0  &  0  &  0      \\
0  &  0  &  J_{w_3}^{-1}   &  0  &  0  &  0      \\
0  &  0  &  0  &  0  &  0  &  0      \\
0  &  0  &  0  &  0  &  0  &  0      \\
0  &  0  &  0  &  0  &  0  &  0      
\end{array} \right] 
\left[ \begin{array}{c}
t_{w_1} \\
t_{w_1} \\
t_{w_1} \\
m_1 \\
m_2 \\
m_3 
\end{array} \right] 
+\left[ \begin{array}{c}
t_{d_1}/J_1 \\
t_{d_2}/J_2 \\
t_{d_3}/J_3 \\
0 \\ 0 \\ 0 \\ 0 \\ 0 \\ 0
\end{array} \right] 
:=\A\x+\B\u+\d.
\label{linearModel1}
\end{eqnarray}
\normalsize

It is worthwhile to notice that (\ref{linearModel1}) is in 
general a time-varying system. The time-variation of the 
system arises from an approximately periodic function of 
$\b(t)=\b(t+P)$, where 
\begin{equation}
P=\frac{2\pi}{\omega_0}=2\pi \sqrt{\frac{a^3}{GM}}
\label{period}
\end{equation}
is the orbital period, $a$ is the orbital radius (approximately 
equal to the spacecraft altitude plus the radius of the Earth), 
and $GM=3.986005*10^{14}{m^3/s^2}$ \cite{wertz78}. This magnetic field $\b(t)$ can be approximately expressed as follows 
\cite{psiaki01}:
\begin{equation}
\left[  \begin{array}{c} 
b_1(t) \\
b_2(t) \\
b_3(t) 
\end{array}  \right] 
= \frac{\mu_f}{a^{3}}
\left[  \begin{array}{c} 
\cos(\omega_0t)\sin(i_m) \\
-\cos(i_m) \\
2\sin(\omega_0t)\sin(i_m) 
\end{array}  \right],
\label{field}
\end{equation}
where $i_m$ is the inclination of the spacecraft orbit with 
respect to the magnetic equator, $\mu_f=7.9 \times 10^{15}$
Wb-m is the field's dipole strength.
The time $t=0$ is measured at the ascending-node
crossing of the magnetic equator. Therefore, the periodic time-varying 
matrix $\B$ in (\ref{linearModel1}) can be written as
\begin{equation}
 \B= \left[ \begin{array}{cccccc}
-J_1^{-1}   &  0  &  0  &  0  & \frac{2\mu_f}{a^3 J_{1}} \sin(i_m) \sin(\omega_0t)  & \frac{\mu_f}{a^3 J_{1}} \cos(i_m)        \\
0  &  -J_2^{-1}   &  0  &   -\frac{2\mu_f}{a^3 J_{2}} \sin(i_m) \sin(\omega_0t)  &  0  &  \frac{\mu_f}{a^3 J_{2}} \sin(i_m) \cos(\omega_0t)     \\
0  &  0  &  -J_3^{-1}   &  -\frac{\mu_f}{a^3 J_{3}} \cos(i_m)   & -\frac{\mu_f}{a^3 J_{3}} \sin(i_m) \cos(\omega_0t)   &  0       \\
J_{w_1}^{-1}   &  0  &  0  &  0  &  0  &  0       \\
0  &  J_{w_2}^{-1}   &  0  &  0  &  0  &  0      \\
0  &  0  &  J_{w_3}^{-1}   &  0  &  0  &  0      \\
0  &  0  &  0  &  0  &  0  &  0      \\
0  &  0  &  0  &  0  &  0  &  0      \\
0  &  0  &  0  &  0  &  0  &  0      
\end{array} \right].
\label{varyB}
\end{equation}
A special case is when $i_m=0$, i.e., the spacecraft orbit is on
the equator plane of the Earth's magnetic field. In this case,
$\b(t)=[0,-\frac{\mu_f}{a^3},0]^{\Tr}$ is a constant vector
and $\B$ is reduced to a constant matrix given as follows:
\begin{equation}
 \B= \left[ \begin{array}{cccccc}
-J_1^{-1}   &  0  &  0  &  0  &  0  & \frac{\mu_f}{a^3 J_{1}}     \\
0  &  -J_2^{-1}   &  0  &  0  &  0  &  0     \\
0  &  0  &  -J_3^{-1}   &  -\frac{\mu_f}{a^3 J_{3}} & 0   &  0       \\
J_{w_1}^{-1}   &  0  &  0  &  0  &  0  &  0       \\
0  &  J_{w_2}^{-1}   &  0  &  0  &  0  &  0      \\
0  &  0  &  J_{w_3}^{-1}   &  0  &  0  &  0      \\
0  &  0  &  0  &  0  &  0  &  0      \\
0  &  0  &  0  &  0  &  0  &  0      \\
0  &  0  &  0  &  0  &  0  &  0      
\end{array} \right].
\label{constantB}
\end{equation}
In the remainder of the discussion, we will consider the discrete
time system of (\ref{linearModel1}) because it 
is more suitable for computer controlled system implementations.
The discrete time system is given as follows:
\begin{eqnarray}
\x_{k+1} = \A \x_k + \B_k \u_k + \d_k.
\label{discreteState}
\end{eqnarray}
Assuming the sampling time is $t_s$, the simplest but less 
accurate discretization formulas to get $\A_k$ and $\B_k$ are
given as follows:
\begin{eqnarray}
\A_k= (\I+ t_s \A), \hspace{0.1in}
\B_k = t_s \B(kt_s).
\label{ABk}
\end{eqnarray}
A slightly more complex but more accurate discretization formulas 
to get $\A_k$ and $\B_k$ are given as follows \cite[page 53]{lvs12}:
\begin{eqnarray}
\A_k= e^{\A t_s}, \hspace{0.1in}
\B_k = \int_0^{t_s} e^{\A \tau} B(\tau) d \tau.
\label{ABk1}
\end{eqnarray}

\section{The LQR design} 

Given the linearized spacecraft model (\ref{linearModel1}) which 
has the state variables composed of spacecraft quaternion $\q$, 
the spacecraft rotational rate with respect to the LVLH frame 
$\boldsymbol{\omega}$, and the reaction wheel rotational speed 
$\boldsymbol{\Omega}$, we can see that to control the spacecraft 
attitude and to manage the reaction wheel momentum
are equivalent to minimize the following objective function
\begin{equation}
\int_0^{\infty} (\x^{\Tr} \Q \x+ \u^{\Tr} \R \u)dt
\label{object}
\end{equation}
under the constraints of (\ref{linearModel1}). This is clearly a 
LQR design problem which has known efficient methods to solve. 
However, in each special case, this system has some
special properties which should be fully utilized to select the most efficient and effective method for each of these cases. 
The corresponding discrete time system is given as follows:
\begin{eqnarray}
\lim_{N \rightarrow \infty} \left( \min  \frac{1}{2} \x_N^{\Tr} \Q_N  \x_N
+\frac{1}{2} \sum_{k=0}^{N-1} \x_k^{\Tr} \Q_k  \x_k+ \u_k^{\Tr} \R_k  \u_k  \right) \nonumber \\
{\mbox s.t.}  \hspace{1in} \x_{k+1} = \A \x_k + \B_k \u_k + \d_k
\label{discreteAll}
\end{eqnarray}

\subsection{Case 1: $\i_m=0$}

It was shown in \cite{yang16} that a spacecraft without a reaction wheel
in this orbit is not controllable. But for a spacecraft with three reaction wheels as we discussed 
in this paper, the system is fully controllable. The controllability condition can be checked straightforward but 
the check is tedious and is omitted in this paper (also the controllability check is not the focus of this paper).
In this case, as we have seen from (\ref{linearModel1}), 
(\ref{constantB}), and (\ref{ABk}) that the linear
system is time-invariant. Therefore, a method for time-varying
system is not appropriate for this simple problem. For 
this linear time-invariant system, the optimal solution of 
(\ref{discreteAll}) is given by (see \cite[page 69]{lvs12})
\begin{equation}
\u_k = -(\R+\B^{\Tr}\P\B)^{-1} \B^{\Tr}  \P \A \x_k = -\K \x_k,
\label{soluC}
\end{equation}
where $\P$ is a constant positive semi-definite solution of 
the following discrete-time algebraic Riccati equation (DARE) 
\begin{equation}
\P =  \Q+\A^{\Tr}\P\A -\A^{\Tr} \P \B(\R+\B^{\Tr}\P\B)^{-1}\B^{\Tr} \P\A.
\label{soluC1}
\end{equation}
There is an efficient algorithms \cite{al84} for this DARE system
and an Matlab function {\tt dare} implements this algorithm.

\subsection{Case 2: $\i_m \ne 0$}

It was shown in \cite{yang16} that a spacecraft without any 
reaction wheel in any orbit of this case is controllable if the 
spacecraft design satisfies some additional conditions imposed 
on $\J$ matrix. By intuition, the system is also controllable
by adding reaction wheels. As a matter of fact, adding reaction 
wheels will achieve better performance of spacecraft attitude 
as pointed in \cite[page 19]{wertz78}. The best algorithm for
this case is a little tricky because $\B$ is a time-varying 
matrix but $\A$ is a constant matrix. Therefore, a method
for time-varying system must be used. The optimal solution 
of (\ref{discreteAll}) is given by (see \cite{bittanti91})
\begin{equation}
\u_k = -(\R_k+\B_k^{\Tr}\P_{k+1}\B_k)^{-1} \B^{\Tr}_k  \P_{k+1} \A_k \x_k = -\K_k \x_k,
\label{soluV}
\end{equation}
where $\P_{k}$ is a periodic positive semi-definite solution of 
the following periodic time-varying Riccati (PTVR) equation
\begin{equation}
\P_k =  \Q_k+\A^{\Tr}_k\P_{k+1}\A_k -\A_k^{\Tr} \P_{k+1} \B_k(\R_k+\B_k^{\Tr}\P_{k+1}\B_k)^{-1}\B_k^{\Tr} \P_{k+1}\A_k.
\label{soluV1}
\end{equation}
Hench and Laub \cite{hl94} developed an efficient algorithm for 
solving the general PTVR equation. However, since $\A_k=\A$ is a 
constant matrix, their algorithm is not optimized. A more 
efficient algorithm in this case was recently proposed in 
\cite{yang16a}, which is particularly useful for time-varying 
system with long period and a constant $\A$ matrix
because it may save hundreds of matrix inverses. The
algorithm is presented below (its proof is in \cite{yang16a}):
\begin{algorithm}  {\ } \\
 \begin{itemize}
\item[] Data: $i_m$, $\J$, $\J_w$, $\Q$, $\R$, the altitude of
the spacecraft (for the calculation of $a$ in (\ref{field})), 
$t_s$ (the selected sample time period), and $p$ 
(the total samples in one period $P=\frac{2\pi}{\omega_0}$).
\item[] Step 1: For $k=1, \ldots, p$, calculate $\A_k$ and $\B_k$ 
using (\ref{ABk}) or (\ref{ABk1}).
\item[] Step 2: Calculate $\E_k$ and $\F_k$ using
\begin{equation}
\E_k =  \left[ \begin{array}{cc} 
\I & \B_k \R^{-1} \B_k^{\Tr}  \\ \0  & \A^{\Tr}
\end{array} \right],
\label{Ek}
\end{equation}
\begin{equation}
\F_k =  \left[ \begin{array}{cc} 
\A & \0  \\ -\Q  & \I
\end{array} \right]=\F.
\label{Fk}
\end{equation}
\item[] Step 3: Calculate $\boldsymbol{\Gamma}_k$, for 
$k=1, \ldots, p$, using
\begin{equation}
\boldsymbol{\Gamma}_k
= \F^{-1} \E_{k} \F^{-1} \E_{k+1} \ldots, 
\F^{-1} \E_{k+p-2} \F^{-1} \E_{k+p-1}.
\label{Gamma}
\end{equation}
\item[] Step 4: Use Schur decomposition
\begin{equation}
\left[ \begin{array}{cc} \W_{11k} &  \W_{12k}  \\
\W_{21k}   &    \W_{22k}
\end{array} \right]^{\Tr} \boldsymbol{\Gamma_k} 
\left[ \begin{array}{cc} \W_{11k} &  \W_{12k}  \\
\W_{21k}   &    \W_{22k}
\end{array} \right]=
\left[ \begin{array}{cc}
\S_{11k} &  \S_{12k}  \\
\0   &    \S_{22k}
\end{array} \right].
\label{decomp1}
\end{equation}
\item[] Step 5: Calculate $\P_k$ using 
\begin{equation}
\P_k =\W_{21k}\W_{11k}^{-1}
\label{solution}
\end{equation}
\end{itemize}
\label{desat}
\end{algorithm} 

\begin{remark}
This algorithm makes full use of the fact that $\A$ is a constant
matrix in (\ref{Fk}). Therefore, $\F$ is a constant matrix and 
the inverse of $\F$ in (\ref{Gamma}) does not need to be 
repeated many times which is the main difference between 
the method in \cite{yang16a} and the method in \cite{hl94} 
(where $\E_k=\E$ is a constant matrix but $\F_k$ is a series of 
time varying matrices and inverse has to take for every 
$\F_k$ with $k=1, \ldots, p$).
\end{remark}

\begin{remark}
The proposed method can easily be extended to the case of using
momentum wheel where the speed of the flywheel is desired
to be a non-zero constant. Let $\boldsymbol{\bar{\Omega}}$ 
be the desired speed of the momentum wheels and 
$\bar{\x}=[ \0^{\Tr},\boldsymbol{\bar{\Omega}}^{\Tr},\0^{\Tr}]^{\Tr}$ . The objective function
of (\ref{object}) should be revised to 
\begin{equation}
\int_0^{\infty} [(\x-\bar{\x})^{\Tr} \Q (\x-\bar{\x})+ \u^{\Tr} \R \u ]dt.
\label{object1}
\end{equation}
\end{remark}

\section{Simulation test} 

Our simulation has several goals. First, we would like to 
show that the proposed design achieves both attitude control and 
reaction wheel momentum management. Second, we would like to 
compare with the design \cite{yang16a} which does not use 
reaction wheels, our purpose is to show that using reaction 
wheels achieves better attitude pointing accuracy. More
important, we would like to demonstrate that the LQR 
design works very well for attitude and desaturation
control for the nonlinear spacecraft in the environment close 
to the reality. Finally, we will discuss the strategy in 
real spacecraft control system implementation.

\begin{figure}[ht]
\centerline{\epsfig{file=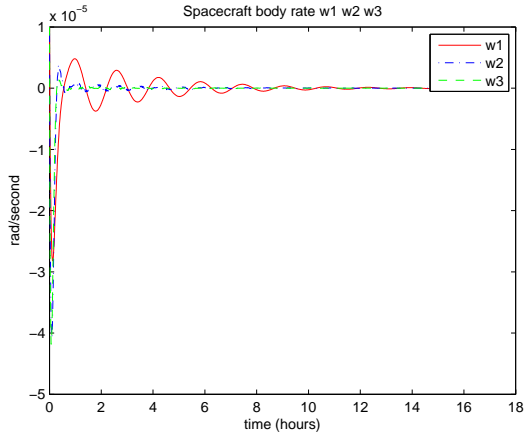,height=6cm,width=8cm}}
\caption{Body rate response $\omega_1$, $\omega_2$, and $\omega_3$.}
\label{fig:linearSCw}
\end{figure}

\begin{figure}[ht]
\centerline{\epsfig{file=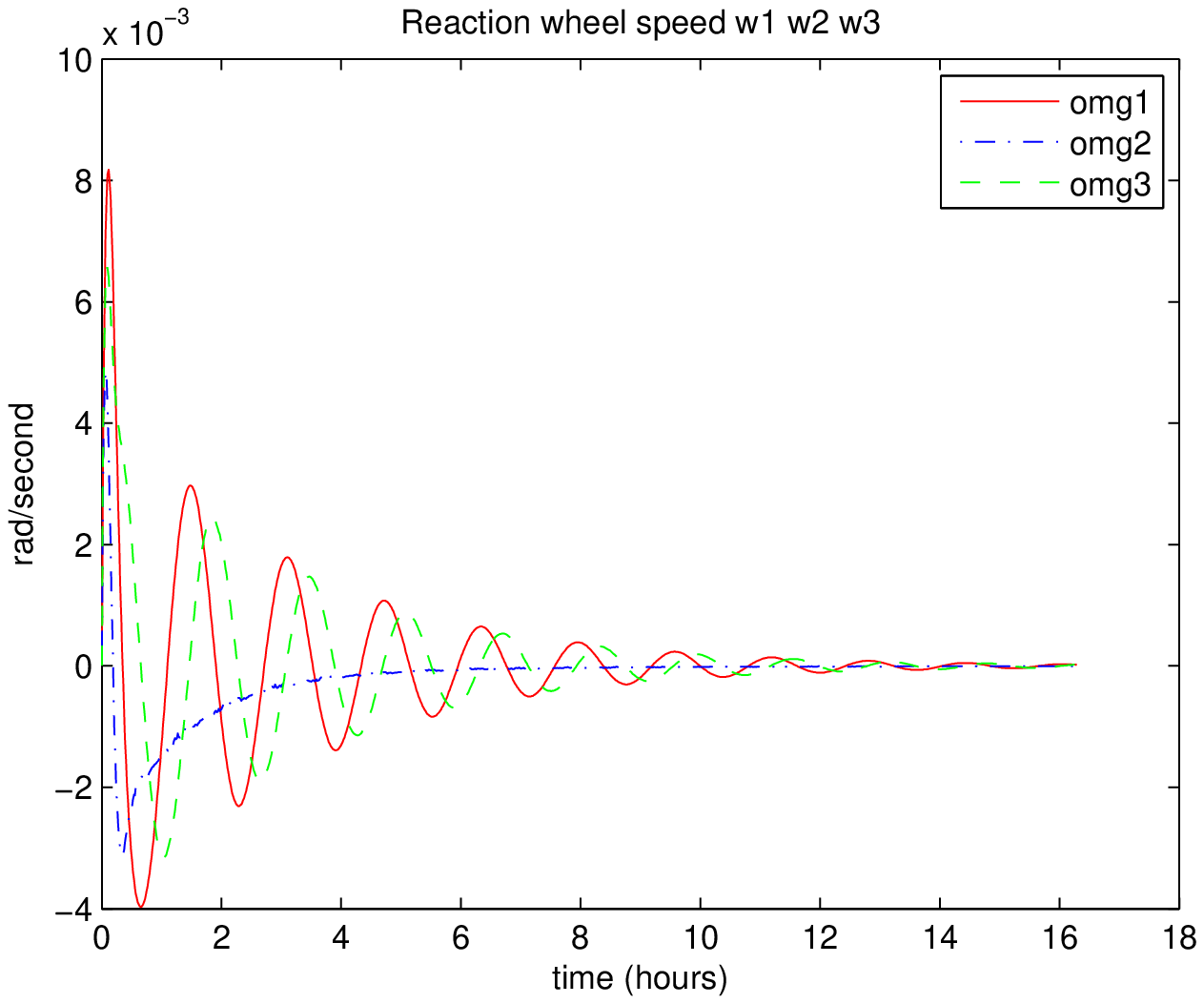,height=6cm,width=8cm}}
\caption{Reaction wheel response $\Omega_1$, $\Omega_2$, and $\Omega_3$.}
\label{fig:linearWheelw}
\end{figure}

\begin{figure}[ht]
\centerline{\epsfig{file=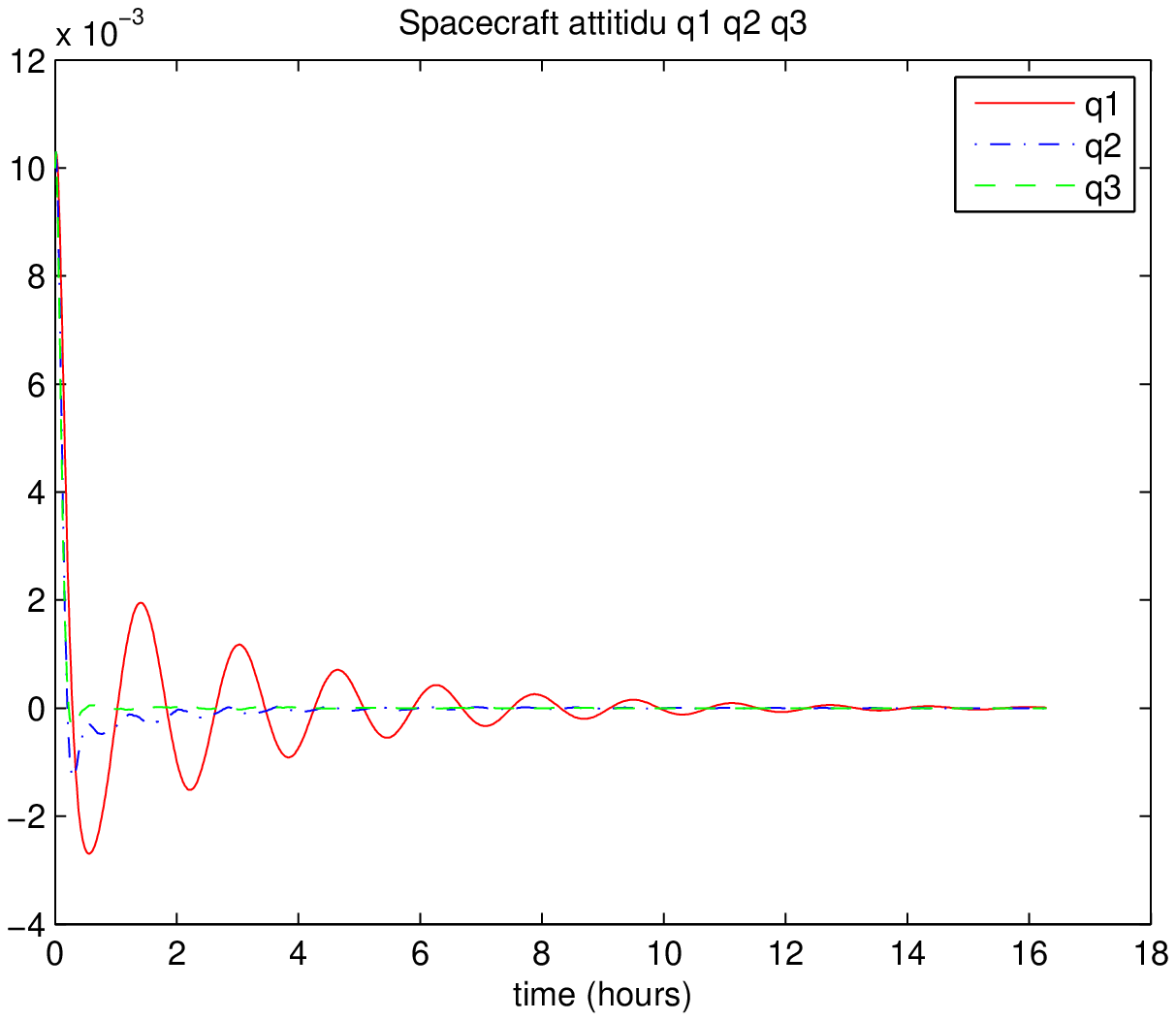,height=6cm,width=8cm}}
\caption{Attitude response $q_1$, $q_2$, and $q_3$.}
\label{fig:linearSCq}
\end{figure}

\begin{figure}[ht]
\centerline{\epsfig{file=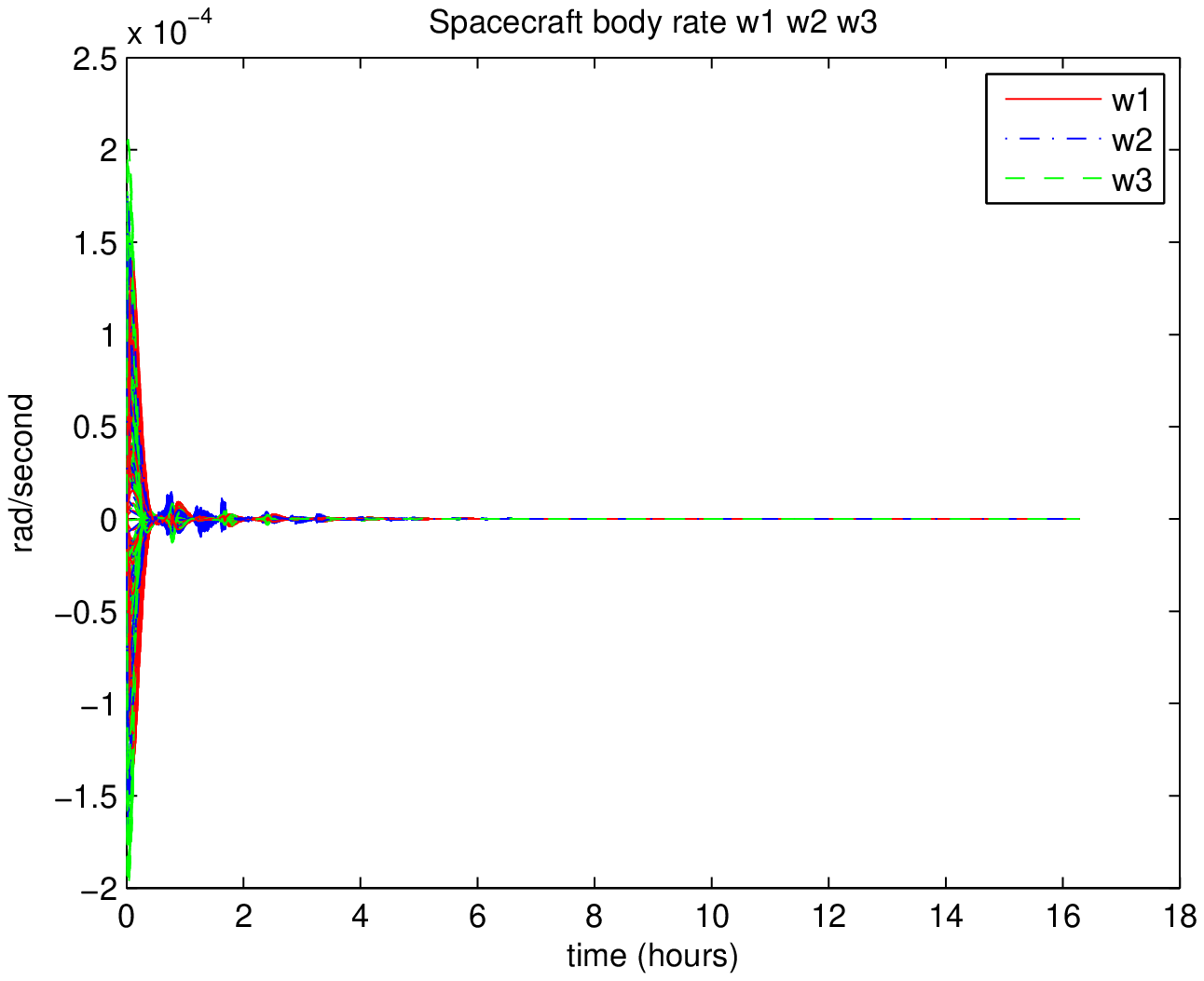,height=6cm,width=8cm}}
\caption{Body rate response $\omega_1$, $\omega_2$, and $\omega_3$.}
\label{fig:NonlinearSCw}
\end{figure}

\begin{figure}[ht]
\centerline{\epsfig{file=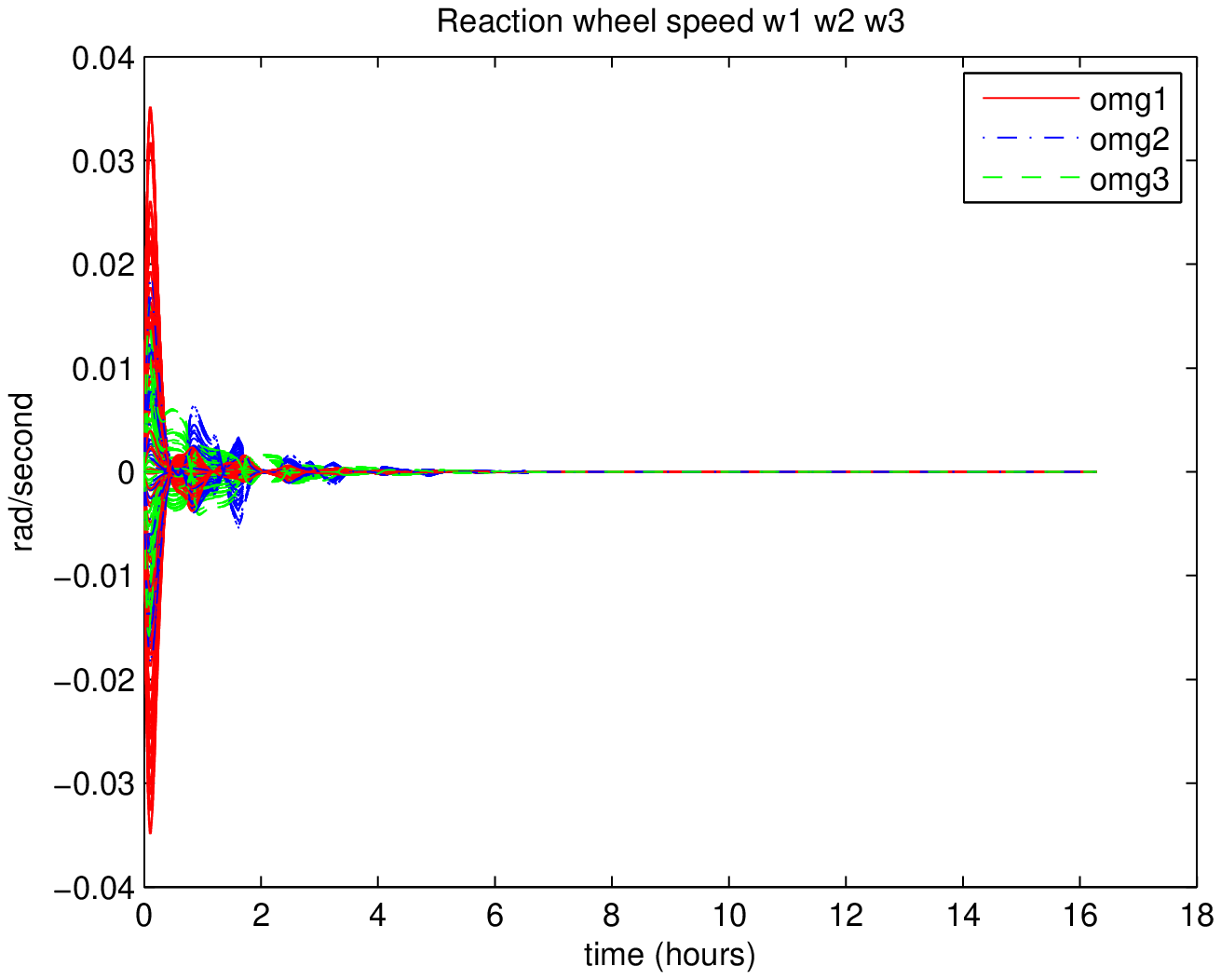,height=6cm,width=8cm}}
\caption{Reaction wheel response $\Omega_1$, $\Omega_2$, and $\Omega_3$.}
\label{fig:NonlinearWheelw}
\end{figure}

\begin{figure}[ht]
\centerline{\epsfig{file=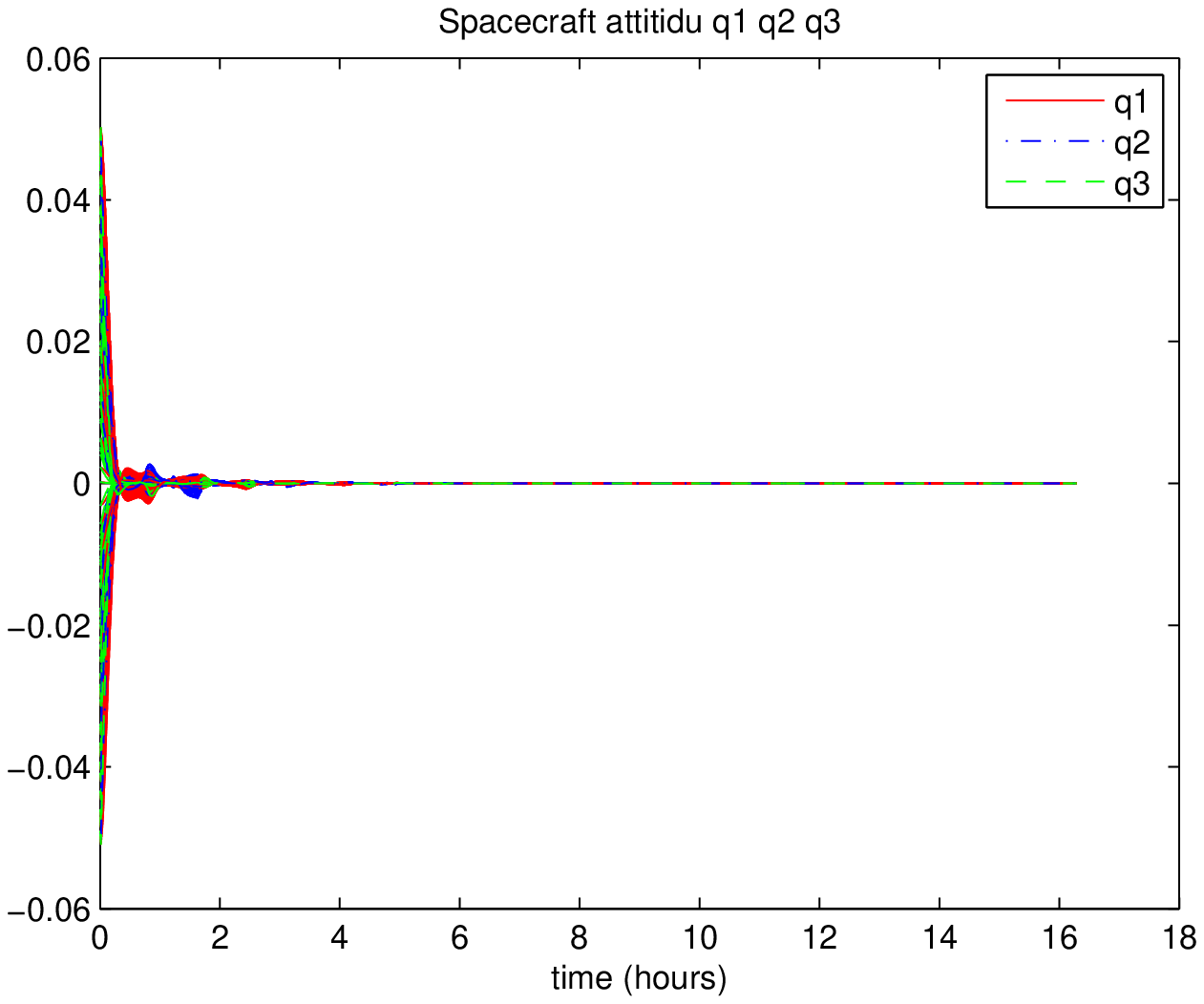,height=6cm,width=8cm}}
\caption{Attitude response $q_1$, $q_2$, and $q_3$.}
\label{fig:NonlinearSCq}
\end{figure}

\subsection{Comparison with the design without reaction wheels}

The proposed design algorithm has been tested using the same 
spacecraft model and orbit parameters as in \cite{yang16a} with 
the spacecraft inertia matrix given by
\[
\J=\diag \left( 250, 150, 100 \right) kg \cdot m^2.
\]
The orbital inclination $i_m=57^o$ and the orbit is assumed to
be circular with the altitude $657$ km. In view of equation 
(\ref{period}), the orbital period is $5863$ seconds and the orbital rate is $\omega_0=0.0011$ rad/second. Assuming that 
the total number of samples taken in one orbit is $100$, then, 
each sample period is $58.6352$ second. It is easy to see that 
all parameters are selected the same as \cite{yang16a} so that 
we can compare the two different designs. Select 
$\Q=\diag([0.001, 0.001, 0.001, 0.001, 0.001, 0.001, 0.02, 0.02, 0.02])$
and $\R=\diag([10^{3}, 10^{3}, 10^{3}, 10^{2}, 10^{2}, 10^{2}])$. 
We have calculated and stored $\P_k$ for $k=0,1,2,\ldots,99$ 
using Algorithm \ref{desat}. Assuming that the initial 
quaternion error is $(0.01, 0.01, 0.01)$, initial body rate 
vector is $(0.00001, 0.00001, 0.00001)$ radians/second, and the
initial wheel speed vector is $(0.00001, 0.00001, 0.00001)$ radians/second,
applying the feedback (\ref{soluV}) to the linearized system 
(\ref{linearModel1}) and (\ref{varyB}),
we get the linearized spacecraft rotational rate response 
described in Figure \ref{fig:linearSCw}, the reaction wheel response descried in 
Figure \ref{fig:linearWheelw}, and the spacecraft attitude responses given in Figures \ref{fig:linearSCq}.

Comparing the response obtained here using both reaction wheels and 
magnetic torque coils and the response obtained in \cite{yang16a} 
that uses magnetic torques only, we can see that both control methods 
stabilize the spacecraft, but using reaction wheels achieve 
much accurate nadir pointing. Also reaction wheel
speeds approach to zero as $t$ goes to infinity. Therefore, the second
design goal for reaction wheel desaturation is achieved nicely.

\subsection{Control of the nonlinear system}

It is nature to ask the following question: can the designed 
controller (\ref{soluV}), which is based on the linearized model, 
stabilize the original nonlinear spacecraft system 
(\ref{dynamicsAll}) with satisfied performance? We answer this 
question by applying the designed controller to the original 
nonlinear spacecraft system (\ref{dynamicsAll}). More specifically,
the LVLH frame rotational rate $\boldsymbol{\omega}_{lvlh}^b$
is calculated using the accurate nonlinear formula (\ref{n9}) 
not the approximated linear model (\ref{9}). The gravity gradient torque $\t_g$ is calculated using the accurate nonlinear formula 
(\ref{ginterM1}) not the approximated linear model (\ref{22}).
The Earth's magnetic field is calculated using the much accurate
International Geomagnetic Reference Field (IGRF) model 
\cite{Finlay15} not the simplified model (\ref{field}). This 
is done as follows. Given the altitude of the spacecraft 
($657$ km), the orbital radius $R$ is $7028$ kilo meters and 
the lateral speed of the spacecraft is $v=R \omega_0$ 
\cite[page 109]{sidi97}.
Assuming that the ascending node at $t=0$ (``now'') is the 
$\X$ axis of the ECEF frame, the velocity vector
$\v = [0,v\cos(i_m),v\sin(i_m)]^{\Tr}$. Using Algorithm 3.4 of
\cite[page 142]{curtis05}, we can get the spacecraft coordinate
in ECI frame at any time after $t=0$. Converting ECI coordinate
to ECEF coordinate, we can calculate a much accurate Earth 
magnetic field vector $\b$ using IGRF model \cite{Finlay15},
which has been implemented in Matlab. Applying this Earth 
magnetic field vector $\b$ and feedback control $\u_k=-\K_k \x_k$ 
designed by the LQR method to (\ref{dynamicsAll}), we control
the nonlinear spacecraft system using the LQR controller.
Also, we allow randomly generated larger initial errors (possibly
$10$ time large than we used in the previous simulation test)
in this simulation test.

The nonlinear spacecraft system response to the LQR controller 
is given in Figures \ref{fig:NonlinearSCw}, 
\ref{fig:NonlinearWheelw}, and \ref{fig:NonlinearSCq}. 
These figures show that the proposed design does 
achieve our design goals. Moreover, the difference between the
linear (approximate) system response and nonlinear (true) system response for the LQR design is very small!

\subsection{Implementation to real system}

In real space environment, even the magnetic field vector 
obtained from the high fidelity IGRF model may not be identical 
to the real magnetic field vector which can be measured by 
magnetometer installed on spacecraft. Therefore, it is suggested
to use the measured magnetic field vector $\b$ to form $\B_k$
in the state feedback (\ref{soluV}). Because of the interaction 
between the magnetic torque coils and the magnetometer, it is a common practice that measurement and control are not taken at 
the same time (some time slot in the sample period is allocated 
to the measurement and the rest time in the sample period is 
allocated for
control). Therefore, a scaling for the  control gain should be 
taken to compensate for the time loss in the sample period
when measurement is taken. For example, if the magnetic field
measurement uses half time of the sample period, the control 
gain should be doubled because only half sample period is used for 
control. This is similar to the method used in \cite{yang14b}.



\section{Conclusions} 
 
In this paper, we developed a reduced quaternion spacecraft model
which includes gravity gradient torque, geomagnetic field along
the spacecrat orbit and its interaction with the magnetic torque 
coils, and the reaction wheels. We investigate a time-varying
LQR design method to control the spacecraft attitude to align 
the body frame with the local vertical local 
horizontal frame and to desaturate the reaction momentum at the
same time. A periodic optimal controller is proposed for this
purpose. The periodic controller design is based on an efficient algorithm
for the periodic time-varying Riccati equations. Simulation test
is given to show that the design objective is achieved and 
the control system using both reaction wheels and magnetic 
torques accomplishes more accurate attitude than the control 
system using only magnetic torques.


\end{document}